\documentclass[a4paper,12pt]{article}
\usepackage{dsfont}
\usepackage{amssymb}
\usepackage{latexsym}
\usepackage{amsmath}
\usepackage{xcolor}
\usepackage{comment}
\usepackage{amsthm}
\usepackage[dvips]{graphicx}
\usepackage{layout} 
\usepackage{ulem}
\usepackage{enumerate}
\usepackage{bm}
\usepackage{bbm} 
\usepackage[bbgreekl]{mathbbol} 
\usepackage{authblk}
\usepackage{setspace} 
\newif\ifrs
\rstrue
\ifrs \usepackage{mathrsfs} \fi  
\newif\ifcol
\coltrue 

\newtheorem{theorem}{Theorem}[section]

\newtheorem{proposition}[theorem]{Proposition}

\newtheorem{remark}[theorem]{Remark}

\numberwithin{equation}{section}
\newtheorem{theorem*}{Theorem}
\newtheorem{ass*}[theorem*]{Assumption}
\newtheorem{note*}[theorem*]{Note}
\newtheorem{lemma*}[theorem*]{Lemma}
\newtheorem{definition*}[theorem*]{Definition}
\newtheorem{proposition*}[theorem*]{Proposition}
\newtheorem{corollary*}[theorem*]{Corollary}
\newtheorem{remark*}[theorem*]{Remark}
\newtheorem{example*}[theorem*]{Example}
\numberwithin{equation}{section}
\newif\ifcol
\colfalse
\ifcol
\newcommand{\colorr}{\color[rgb]{0.8,0,0}}

\newcommand{\colorb}{\color[rgb]{0,0,0.8}}

\newcommand{\colorn}{\color[rgb]{1,1,1}}

\newcommand{\coloroy}{\color[rgb]{1,0.95,0}}

\else
\newcommand{\colorb}{\color{black}}
\newcommand{\colorr}{\color{black}}

\newcommand{\colorn}{\color{black}}

\newcommand{\coloroy}{\color{black}}
\fi
\newif\ifcol
\colfalse
\ifcol
\newcommand{\sred}{\color[rgb]{0.8,0,0}}

\newcommand{\sblue}{\color[rgb]{0,0,0.8}}
\else
\newcommand{\sred}{\color{black}}
\newcommand{\sblue}{\color{black}}
\fi
\newif\ifcol
\colfalse
\ifcol

\else
\fi
\excludecomment{en-text}
\includecomment{jp-text}
\includecomment{comment}
\def\ol{\overline}

\def\infm{{\infty\text{--}}}
\def\inftym{\infm}

\def\bd{\begin{description}}
\def\ed{\end{description}}

\def\D2{\bbD_{2,\infty-}}

\def\D{{\bf D}}

\def\calf{{\cal F}}
\def\calg{{\cal G}}

\def\ds{\displaystyle}
\def\yeq{\>=\>}
\def\yleq{\>\leq\>}
\def\ygeq{\>\geq\>}

\def\sfp{{\sf p}}

\def\simleq{\ \raisebox{-.7ex}{$\stackrel{{\textstyle <}}{\sim}$}\ }

\def\ep{\epsilon}
\def\half{\frac{1}{2}}

\def\down{\downarrow}

\def\halflineskip{\vspace*{3mm}}
\def\nn{\nonumber}
\def\be{\begin{equation}}
\def\ee{\end{equation}}
\def\bea{\begin{eqnarray}}
\def\eea{\end{eqnarray}}
\def\beas{\begin{eqnarray*}}
\def\eeas{\end{eqnarray*}}
\def\bi{\begin{itemize}}
\def\ei{\end{itemize}}
\def\im{\item}
\def\bd{\begin{description}}
\def\ed{\end{description}}
\def\dotc{\stackrel{\circ}{C}}

\newcommand{\bbB}{{\mathbb B}}

\newcommand{\bbD}{{\mathbb D}}
\newcommand{\bbE}{{\mathbb E}}

\newcommand{\bbH}{{\mathbb H}}

\newcommand{\bbN}{{\mathbb N}}

\newcommand{\bbP}{{\mathbb P}}

\newcommand{\bbR}{{\mathbb R}}

\newcommand{\bbT}{{\mathbb T}}
\newcommand{\bbU}{{\mathbb U}}
\newcommand{\bbV}{{\mathbb V}}

\newcommand{\bbY}{{\mathbb Y}}
\newcommand{\bbZ}{{\mathbb Z}}

\begin{document}

\title{
Simplified quasi-likelihood analysis for a locally asymptotically quadratic random field
\footnote{
This work was in part supported by 
Japan Science and Technology Agency CREST JPMJCR14D7, JPMJCR2115; 
Japan Society for the Promotion of Science Grants-in-Aid for Scientific Research 
No. 17H01702 (Scientific Research);  
and by a Cooperative Research Program of the Institute of Statistical Mathematics. 
The author thanks Junichiro Yoshida for valuable comments. 
}
}
\author{Nakahiro Yoshida}
\affil{Graduate School of Mathematical Sciences, University of Tokyo
\footnote{Graduate School of Mathematical Sciences, University of Tokyo: 3-8-1 Komaba, Meguro-ku, Tokyo 153-8914, Japan. e-mail: nakahiro@ms.u-tokyo.ac.jp}\\
Japan Science and Technology Agency CREST\\
The Institute of Statistical Mathematics}
\date{
December 27, 2021
}
\maketitle
\noindent
{\it Summary} 
The asymptotic decision theory by Le Cam and H\'ajek has been given a lucid perspective by the Ibragimov-Hasminskii theory on convergence of the likelihood random field. 
Their scheme has been applied to stochastic processes by Kutoyants, and today this plot is called the IHK program. 
This scheme ensures that asymptotic properties of an estimator follow directly from the convergence of the random field if a large deviation estimate exists. 
The quasi-likelihood analysis (QLA) proved a polynomial type large deviation (PLD) inequality to go through a bottleneck of the program. 
A conclusion of the QLA is that if the quasi-likelihood random field is asymptotically quadratic and if a key index reflecting identifiability the random field has is non-degenerate, then the  PLD inequality is always valid, and as a result, the IHK program can run. 
Many studies already took advantage of the QLA theory. 
However, not a few of them are using it in an inefficient way yet. 
The aim of this paper is to provide a reformed and simplified version of the QLA and 
to improve accessibility to the theory. 
As an example of the effects of the program and 
the PLD, the user can obtain asymptotic properties of the quasi-Bayesian estimator by only verifying non-degeneracy of the key index. 
\ \\
\ \\
{\it Keywords and phrases } 
Ibragimov-Has'minskii theory, 
quasi-likelihood analysis, 
polynomial type large deviation, 
random field, 
asymptotic decision theory, non-ergodic statistics. 


\section{Introduction}
The asymptotic decision theory by Le Cam and H\'ajek has been given a lucid perspective by the Ibragimov-Has'minskii theory 
(\cite{ibragimov1973asymptotic,ibragimov1973asymptoticbayes,IbragimovHascprimeminskiui1981}) on convergence of the likelihood random field. 
Their scheme has been applied to stochastic processes by Kutoyants 
(\cite{Kutoyants1984,Kutoyants1994,Kutoyants1998,Kutoyants2004}), and today this plot is called the IHK program. 
This scheme ensures that asymptotic properties of an estimator follow directly from the convergence  of the random field if a large deviation estimate exists. 

The theory of statistical inference for stochastic processes is heading toward 
heavily dependent stochastic systems: 
nonlinear, non-Markovian, non-stationary, non-ergodic, discrete and/or dependent sampling. 
A formal extension of the likelihood analysis to a quasi-likelihood analysis is 
inevitable. However, the large deviation inequality is an issue even within the likelihood analysis for stochastic processes. 
The quasi-likelihood analysis (QLA) proved 
a polynomial type large deviation (PLD) inequality to go through a bottleneck of the program 
(Yoshida \cite{Yoshida2011}).
\footnote{The term ``quasi-likelihood'' is not in the sense of GLM.
We use ''quasi-likelihood analysis'' because statistical inference for sampled stochastic processes cannot avoid a quasi-likelihood function for estimation. 
The method is relatively new, but not because of ``quasi''. The difficulty in large deviation estimates already existed in the likelihood analysis for stochastic processes.  
}
As a conclusion of the QLA theory, if the quasi-likelihood random field is locally asymptotically quadratic (LAQ) and if a key index reflecting identifiability the random field has is non-degenerate, then the polynomial type large deviation inequality is always valid, and as a result, the IHK program can run.   

Since an ad hoc model-dependent method is not necessary, 
the QLA is universal and can apply to various dependent models. 
Many studies are based on and taking advantage of the QLA. 
These applications include 
sampled ergodic diffusion processes (Yoshida \cite{Yoshida2011}), 
adaptive estimation for diffusion processes (Uchida and Yoshida \cite{UchidaYoshida2012Adaptive}), 
adaptive Bayes type estimators for ergodic diffusion processes (Uchida and Yoshida \cite{uchida2014adaptive}), 
approximate self-weighted LAD estimation of discretely observed ergodic Ornstein-Uhlenbeck processes (Masuda \cite{Masuda2010a}), 
parametric estimation of L{\'e}vy processes (Masuda \cite{masuda2015parametric}),
Gaussian quasi-likelihood random fields for ergodic L\'evy driven SDE (Masuda \cite{masuda2013convergence}), 
and 
ergodic point processes for limit order book (Clinet and Yoshida \cite{clinet2015statistical}). 
%
Thanks to its flexibility, the QLA is also applicable to non-ergodic statistics: 
volatility parameter estimation in regular sampling of finite time horizon (Uchida and Yoshida \cite{UchidaYoshida2013})  
and in non-synchronous sampling (Ogihara and Yoshida \cite{ogihara2014quasi}), 
a non-ergodic point process regression model (Ogihara and Yoshida \cite{ogihara2015quasi}). 
Analysis of complex algorithms is possible by relying on  
the universal design of the QLA: 
hybrid multi-step estimators (Kamatani and Uchida \cite{KamataniUchida2014}), 
adaptive Bayes estimators and hybrid estimators for small diffusion processes based on sampled data (Nomura and Uchida \cite{nomura2016adaptive}). 
%
%
Information criteria, sparse estimation and regularization methods are recently understood 
in the framework of the QLA: 
contrast-based information criterion for diffusion processes (Uchida \cite{Uchida2010}), 
AIC for non-concave penalized likelihood method (Umezu et al. \cite{umezu2019aic}), 
Schwarz type model comparison for LAQ models (Eguchi and Masuda \cite{eguchi2016schwarz}), 
moment convergence of regularized least-squares estimator for linear regression model (Shimizu \cite{shimizu2017moment}),  
moment convergence in regularized estimation under multiple and mixed-rates asymptotics (Masuda and Shimizu \cite{masuda2017moment}), 
penalized method and polynomial type large deviation inequality (Kinoshita and Yoshida \cite{kinoshita2019penalized}) 
and related 
Suzuki and Yoshida (\cite{suzuki2020penalized}). 
Jump filtering problems: 
jump diffusion processes Ogihara and Yoshida(\cite{OgiharaYoshida2011}), 
threshold estimation for stochastic processes with small noise (Shimizu \cite{shimizu2017threshold}), 
global jump filters (Inatsugu and Yoshida \cite{inatsugu2021global}). 
Partial quasi-likelihood analysis: 
Yoshida \cite{yoshida2018partial}. 
Such variety of applications are demonstrating the universality of the framework of the QLA. 
Since the IHK program runs there, we can obtain 
limit theorems and the $L^p$-boundedness of the QL estimators (quasi-maximum likelihood estimator and the quasi-Bayesian estimator), 
which is indispensable to develop statistical theories. 

The essence of the QLA is the polynomial type large deviation inequality 
that was proved in a general setting  
(Yoshida \cite{Yoshida2011}). 
\begin{en-text}
We say the quasi-log likelihood random field is 
locally asymptotically quadratic (LAQ) if 
it is approximated by 
a random quadratic function in the local parameter. 
\end{en-text}
Since the LAQ property quite often appears when the model is differentiable, 
Yoshida \cite{Yoshida2011} was based on this structure. 
Because of it, the limit distribution of the associated estimators 
has an explicit expression. 
The paper \cite{Yoshida2011} gave it, but due to a general way of writing, 
not a few users are apt to avoid following that passage after the PLD's theorem 
and try to reconstruct it in each situation. 
However, such a task is unnecessary in fact. 
%
%
Besides, four time differentiability is often assumed in many applications of the QLA.  
It may be only because a handy condition in \cite{Yoshida2011} assumed 
an estimate of the supremum of the third-order derivative of the quasi-log likelihood random field 
$\bbH_T$, though the paper gave a condition ($[A1']$) to treat $\bbH_T$ of class $C^2$. 

The aim of this paper is to provide a simplified version of the QLA theory 
directly connecting the assumptions with the limit theorems 
in order to improve accessibility to the theory. 
%
\begin{en-text}
Verifying non-degeneracy of a key index is essentially the only task requested to the user 
and in particular it becomes trivial in ergodic statistics. 
\end{en-text}
Essentially, the user is only requested to verify non-degeneracy of a key index, and this task is trivial in particular in ergodic statistics. 
%
We will give handy conditions for the quasi-likelihood random field 
of class $C^2$, based on \cite{Yoshida2011}, in order to reach the asymptotic properties of the estimators at a single leap. 
%
Some assumptions in \cite{Yoshida2011} are arranged and replaced by simple-looking ones in this paper. 
This simplification will serve for future progress e.g. in analysis of regularization methods.  
The LAQ property we adopted here is just one principle of separation, and 
it is possible to develop a similar theory for a non-LAQ type random field; 
see Kinoshita and Yoshida \cite{kinoshita2019penalized} for a case of regularization. 
%

A smart way of presenting the theory is to use the convergence of 
the quasi-likelihood random field $\bbZ_T$ to a random field $\bbZ$ 
in the function space $\widehat{C}(\bbR^\sfp)$, the separable Banach space of 
continuous functions $f$ on $\bbR^\sfp$ satisfying $\lim_{|u|\to\infty}f(u)=0$, 
equipped with the supremum norm. 
This plot is possible but to carry out it,  one needs a suitable measurable extension of $\bbZ_T$ to the outside of the originally given local parameter space and an  argument about tightness of random fields on the non-compact $\bbR^\sfp$. 
In this article, we dared avoid this approach to give priority to simplicity. 
As a result, the presentation of the theory is now much more elementary than 
Yoshida \cite{Yoshida2011}. 
{\sred 
Additionally, though they are classical, 
some basic properties like the first order efficiency and asymptotic equivalence between the maximum likelihood estimator 
and the Bayesian estimator will be given within the QLA framework.
For convenience of use, we will give and detail several versions of theorems and conditions 
in different situations. 
}


\begin{en-text}
The QLA is a framework of statistical inference for stochastic processes. 
It features the polynomial type large deviation of the quasi likelihood random field. 
Through QLA, one can systematically derive 
limit theorems and precise tail probability estimates of the associated QLA estimators such as quasi maximum likelihood estimator (QMLE), 
quasi Bayesian estimator (QBE) and various adaptive estimators. 
The importance of such precise estimates of tail probability is well recognized in asymptotic decision theory, 
prediction, theory of information criteria for model selection, asymptotic expansion, etc. 
\end{en-text}

\section{A simplified QLA in non-ergodic statistics
}
Given a probability space $(\Omega,\calf,P)$ 
and a bounded open set $\Theta$ in $\bbR^\sfp$, 
we consider a random field $\bbH_T:\Omega\times\ol{\Theta}\to\bbR$, 
a function measurable with respect to the product $\sigma$-field 
$\calf\times\bbB(\Theta)$, $\bbB(\Theta)$ being the Borel $\sigma$-field of $\Theta$.
\footnote{Because of the assumptions below about the continuity of $\bbH_T$ and the separability of $\Theta$, this is equivalent to that the function $\bbH_T(\cdot,\theta)$ is measurable for each $\theta\in\Theta$.}
Here $t\in\bbT$, a subset of $\bbR_+=[0,\infty)$ satisfying $\sup\bbT=\infty$. 
We suppose that $\bbH_T$ is continuous and of class $C^2$, that is, 
for every $\omega\in\Omega$, 
the mapping $\Theta\ni\theta\mapsto\bbH_T(\theta)\in\bbR$ is of class $C^2$ 
and that $\bbH_T$ is continuously extended to $\partial\Theta$. 
We shall present a simplified version of 
the polynomial type large deviation inequality of Yoshida \cite{Yoshida2011} 
under a handy set of 
sufficient conditions.

Let $\theta^*\in\Theta$. 
Define $\Delta_T$ and $\Gamma_T(\theta)$ by 
\bea
\Delta_T \yeq \partial_\theta\bbH_T(\theta^*)a_T
&\text{and}&
\Gamma_T(\theta) \yeq -a_T^\star\partial_\theta^2\bbH_T(\theta)a_T
\eea
respectively, where $\star$ denotes the matrix transpose. 
Let $a_T\in GL(\bbR^\sfp)$ be a scaling matrix such that $|a_T|\to0$ as $n\to\infty$. 
We suppose that $\Gamma$ is a $\sfp\times\sfp$ symmetric random matrix. 
Let $U(\theta,r)=\{\theta'\in\bbR^\sfp;\>|\theta'-\theta|<r\}$ for $\theta\in\Theta$ and $r>0$. 
There exists a positive constant ${\tt r}_0$ such that $U(\theta^*,{\tt r}_0)\subset\Theta$. 
The minimum and maximum eigenvalues of the symmetric matrix $M$ are denoted by 
$\lambda_{\text{min}}(M)$ and $\lambda_{\text{max}}(M)$, respectively. 
Let $b_T=\big\{\lambda_{\text{min}}(a_T^\star a_T)\}^{-1}$. 
In particular, $b_T\to\infty$ as $n\to\infty$. 
Moreover, we assume that 
\bea\label{0302220841} 
b_T^{-1}\yleq \lambda_{\text{max}}(a_T^\star a_T) \yleq C_0b_T^{-1}\quad(T\in\bbT)
\eea
for some constant $C_0\in[1,\infty)$. 
A typical case is $n$ for $b_T$, and $n^{-1/2}I_\sfp$ for $a_T$, where $I_\sfp$ is the identity matrix.
\begin{remark}\rm
In an ergodic diffusion model, the parameter $\theta_1$ of the diffusion coefficient and 
the parameter $\theta_2$ of the drift coefficient have different convergence rates 
in estimation with high frequency data. 
Then Condition (\ref{0302220841}) may seem restrictive, but 
it is incorrect. The random field $\bbH_n$ is not necessarily the same as 
a quasi-log likelihood function $\Psi_n$ used for estimation in reality, where $n=T\in\bbT=\bbN$. 
The random field $\bbH_n$ is rather ''living in the proof'' in various manners. 
Consider a joint quasi-maximum likelihood estimator $(\hat{\theta}_{1,n},\hat{\theta}_{2,n})$ for $(\theta_1,\theta_2)$. 
To analyze the asymptotic behavior of $\hat{\theta}_{1,n}$, 
the random field 
$\bbH_n(\theta_1)=\Psi_n(\theta_1,\hat{\theta}_{2,n})$ can be used. 
$\bbH_n(\theta_1)$ is estimated by taking supremum about the second argument of 
$\Psi_n$ at some stage. 
For $\hat{\theta}_{2,n}$, one can switch $\bbH_n$ to a different random field 
$\bbH_n(\theta_2)=\Psi_n(\hat{\theta}_{1,n},\theta_2)$. 
Such a stepwise application of the QLA in the present article's form 
can be observed in many studies; 
see Yoshida \cite{Yoshida2011}, Uchida and Yoshida \cite{uchida2012adaptive,uchida2014adaptive} and 
the papers listed in Introduction. 
\end{remark}

Define $\bbY_T:\Omega\times\Theta\to\bbR$ by 
\beas 
\bbY_T(\theta) 
&=& 
\frac{1}{b_T}\big\{\bbH_T(\theta)-\bbH_T(\theta^*)\big\}
\quad(\theta\in\Theta)
\eeas
for {\sred $T\in\bbT$}. 
Let $\bbY:\Omega\times\Theta\to\bbR$ be a continuous random field. 
Let $L>0$ be a positive number. 
\bd
\im[{\bf [S1]}] 
Parameters $\alpha$, $\beta_1$, $\beta_2$, 
$\rho_1$ and $\rho_2$ satisfy 
the following inequalities: 
\beas &&
0<\beta_1<1/2,\quad 
0<\rho_1<\min\big\{1,\alpha/(1-\alpha),2\beta_1/(1-\alpha)\big\},\quad
\\&&
{\colorr 0<\>}2\alpha<\rho_2,\quad
\beta_2\geq0,\quad
1-2\beta_2-\rho_2>0.
\eeas
\ed

\bd
\im[{\bf [S2]}] 
\bd
\im[{\bf (i)}] 
There exists a positive random variable $\chi_0$ and the following conditions are fulfilled. 
\bd
\im[{\bf (i-1)}] 
$
\bbY(\theta) \yeq \bbY(\theta) -\bbY(\theta^*)
\yleq -\chi_0\big|\theta-\theta^*\big|^2
$
for all $\theta\in\Theta$. 
\im[{\bf (i-2)}] 
For some constant $C_L$, it holds that 
\beas 
P\big[\chi_0\leq r^{-(\rho_2-2\alpha)}\big]
&\leq& \frac{C_L}{r^L}\quad(r>0).
\eeas
\ed
\im[{\bf (ii)}] 
For some constant $C_L$, it holds that 
\beas 
P\big[\lambda_{min}(\Gamma)<
r^{-\rho_1}\big]\leq\frac{C_L}{r^L}\quad(r>0)
\eeas
\ed

\ed


Let $\beta=\alpha/(1-\alpha)$. 
Let $\|V\|_p=\big(E[|V|^p])^{1/p}$ for $p>0$ and a matrix valued random variable $V$. 
\bd
\im[{\bf [S3]}]
\bd
\im[{\bf (i)}] For $M_1=L(1-\rho_1)^{-1}$, 
$
\sup_{T\in\bbT}\big\||\Delta_T|\big\|_{M_1}<\infty.
$
\im[{\bf (ii)}] For $M_2=L(1-2\beta_2-\rho_2)^{-1}$, 
{\colorr
\beas 
\sup_{T\in\bbT}
\bigg\|\sup_{\theta\in\Theta\setminus U(\theta^*,b_T^{-\alpha/2})}
b_T^{\half-\beta_2}\big|\bbY_T(\theta)-\bbY(\theta)\big|\bigg\|_{M_2}
&<&\infty.
\eeas
}
\begin{en-text}
\im[{\bf (iii)}] For $M_3=L(\beta-\rho_1)^{-1}$, 
\beas 
\sup_{T\in\bbT}\bigg\|
b_T^{-1}\sup_{\theta\in\Theta}\big|\partial_\theta^3\bbH_T(\theta)\big|
\bigg\|_{M_3}
&<& \infty. 
\eeas
\end{en-text}
{\colorr
\im[(iii)] For $M_3=L(\beta-\rho_1)^{-1}$
\beas 
\sup_{T\in\bbT}\bigg\|\sup_{u\in{\sred\delta^{-1}(\Theta-\theta^*),}\>|u|\leq1}
\big|\Gamma_T(\theta^*+\delta u)-\Gamma_T(\theta^*)\big|\bigg\|_{M_3} &=& O(\delta)\quad(\delta\down0). 
\eeas
}
\im[{\bf (iv)}] For $M_4=L\big(2\beta_1(1-\alpha)^{-1}-\rho_1\big)^{-1}$, 
\beas 
\sup_{T\in\bbT}\big\|
b_T^{\beta_1}\big|\Gamma_T{\colorr(\theta^*)}-\Gamma\big|\big\|_{M_4}
&<& \infty.
\eeas
\ed
\ed

\begin{remark}\label{0302240406}\rm 
(i) 
In the above conditions, each constant $C_L$ is independent of $r$ and $T$, 
but may depend on the parameters appearing in $[S1]$ as well as $L$. 
(ii) 
\begin{en-text}
The condition of three time differentiability of $\bbH_T$ 
is for simplicity of the presentation. 
It is possible to reduce it to the two time differentiability. 
See Yoshida \cite{Yoshida2011}. 
Remark that  
estimation of the supremum of $\partial_\theta^i\bbH_T(\theta)$ 
requires one more differentiability 
if one uses its derivative. 
The use of the fourth derivative is simplifying the argument and maybe sufficient 
in practice, but not the best way theoretically. 
(iii) 
\end{en-text}
In applications, we often need to estimate 
the supremum of a sequence of martingales depending on $\theta$ 
to verify the above moment conditions. 
Use of Sobolev's embedding inequality 
is a simple solution. 
{\sred The GRR inequality is an alternative if one wants to reduce differentiability assumptions.}
(iii) The random matrix $\Gamma$ is positive-definite a.s. 
if $[S2]$ (ii) is satisfied. 
\end{remark}

Let $\bbU_T=\{u\in\bbR^\sfp;\>\theta^*+a_Tu\in\Theta\}$ and $\bbV_T(r)=\{u\in\bbU_T;\>|u|\geq r\}$ for $r>0$. 
Define the random field $\bbZ_T$ on $\bbU_T$ by 
\beas 
\bbZ_T(u) 
&=& 
\exp\big(\bbH_T(\theta^*+a_Tu)-\bbH_T(\theta^*)\big)
\eeas
for $u\in\bbU_T$. 
Following Yoshida \cite{Yoshida2011}, 
we give a polynomial type large deviation inequality for 
the random field $\bbZ_T$. 
\begin{theorem}\label{0302201901}
Given a positive constant $L$, suppose that $[S1]$, $[S2]$ and $[S3]$ are fulfilled. 
Then there exists a constant $C_L$ such that 
\bea\label{202005311807}
P\bigg[\sup_{u\in\bbV_T(r)}\bbZ_T(u)\geq \exp\big(-2^{-1} r^{2-(\rho_1\vee\rho_2)}\big)
\bigg]
&\leq&
\frac{C_L}{r^L}
\eea
for all $r>0$ and $T\in\bbT$. 
The supremum of the empty set should read $-\infty$. 
\end{theorem}
\proof 
Suppose that the constants $\alpha,\beta_1,\beta_2,\rho_1,\rho_2$ satisfy 
Condition $[S1]$. 
We will apply Theorem 1 of \cite{Yoshida2011} with $\rho=2$ 
for $\bbH_T$ of class $C^2$.
According to Section 3.1 of \cite{Yoshida2011}, 
it suffices to verify Conditions $[A1']$ and  $[A2]$-$[A6]$ therein. 
Condition $[A4]$ of \cite{Yoshida2011} with $\rho=2$ and 
the condition that $\alpha\in(0,1)$ 
are satisfied 
under $[S1]$ since $\rho_2<1$. 

Condition $[A1']$ of \cite{Yoshida2011} requires the estimate 
{\sblue
\bea\label{202005311744}
1_{\{r\leq b_T^{(1-\alpha)/2}\}}P\big[S_T'(r)^c\big] &\leq & \frac{C_L}{r^L}\quad(r>0,\>T\in\bbT)
\eea
}
for some constant $C_L$, where 
the event $S_T'(r)$ is defined by 
{\sblue
\bea\label{0312171307} 
S_T'(r) 
&=& 
\left\{\sup_{h:\>\theta^*+h\in\Theta,\atop\>|h|\leq C_0^{1/2}b_T^{-\alpha/2}}
\big|\Gamma_T(\theta^*+h)-\Gamma\big|<r^{-\rho_1}
\right\};
\eea
see Remark \ref{0312300647}. 
}\noindent 
%
To verify (\ref{202005311744}), we may 
{\sblue consider the case where $r\leq b_T^{(1-\alpha)/2}$}, equivalently, 
\bea\label{202005311747}
b_T^{-1} &\leq& {\sblue r^{-2/(1-\alpha)}.} 
\eea 
\begin{en-text}
Otherwise, $S_T'(r)^c=\emptyset$, and there is nothing to show. 
Furthermore, 
we may assume $r$ is sufficiently large (in particular, $r\geq1$) to show Inequality (\ref{202005311744}), by changing $C_L$ if necessary. 
\end{en-text}
%
We have 
\bea\label{202005311841}
{\sblue 1_{\{r\leq b_T^{(1-\alpha)/2}\}}}
P\big[S_T'(r)^c\big]
&\leq&
\bbP_1(T,r)+\bbP_2(T,r),
\eea
where 
\beas
\bbP_1(T,r)
&=&
P\left[
\sup_{h{\sred\in\Theta-\theta^*}:
{\sblue |h|\leq C_0^{1/2}r^{-\beta}}}
\big|\Gamma_T(\theta^*+h)-\Gamma_T(\theta^*)\big|\ygeq\half\> r^{-\rho_1}
\right]
\eeas
and 
\beas
\bbP_2(T,r)
&=&
P\left[\big|\Gamma_T(\theta^*)-\Gamma\big|\ygeq\half\>r^{-\rho_1}\right]
1_{\left\{b_T^{-1} \leq {\sblue r^{-2/(1-\alpha)}}\right\}}
\eeas
in view of (\ref{202005311747}). 
For sufficiently large $r$, by Condition $[S3]$ (iii), 
\bea\label{202005311842}
\sup_{T\in\bbT}\bbP_1(T,r)
&\simleq&
r^{M_3\rho_1}
\sup_{T\in\bbT}E\left[
\sup_{h{\sred\in\Theta-\theta^*}:|h|\leq {\sblue C_0^{1/2}}r^{-\beta}}
\big|\Gamma_T(\theta^*+h)-\Gamma_T(\theta^*)\big|^{M_3}\right]
\nn\\&\simleq&
r^{-M_3(\beta-\rho_1)}
\nn\\&=&
r^{-L}.
\eea
Next, by Condition $[S3]$ (iv), we have 
\bea\label{202005311843}
\sup_{T\in\bbT}\bbP_2(T,r)
&\simleq&
\sup_{T\in\bbT}\bigg(b_T^{-M_4\beta_1}r^{M_4\rho_1}
1_{\left\{b_T^{-1} \leq {\sblue r^{-2/(1-\alpha)}}\right\}}\bigg)
\nn\\&\simleq&
r^{-M_4(2\beta_1/(1-\alpha)-\rho_1)}
\nn\\&=&
r^{-L}.
\eea
%
From (\ref{202005311841}), (\ref{202005311842}) and (\ref{202005311843}), we obtain (\ref{202005311744}), 
therefore 
$[A1']$ of \cite{Yoshida2011} has been verified. 

Condition $[A6]$ of \cite{Yoshida2011} follows from Condition $[S3]$ (i) and 
Condition $[S3]$ (ii). 
Condition $[S2]$ (i) ensures Conditions $[A3]$ for $\rho=2$ 
and $[A5]$ of \cite{Yoshida2011}. 
Moreover, $[S2]$ (ii) verifies $[A2]$ of \cite{Yoshida2011}. 
Now, as already mentioned, we apply Theorem 1 of \cite{Yoshida2011} 
to $\bbZ_T$ for $\rho=2$ in order to obtain (\ref{202005311807}). 
\qed\halflineskip

{\sblue 
\begin{remark}\label{0312300647}\rm
In the proof of Theorem \ref{202005311807}, 
we changed $S_T'(r;\xi_0)$ of Yoshida \cite{Yoshida2011} to (\ref{0312171307}), 
and Condition $[A1']$ therein to (\ref{202005311744}). 
The latter apparently weakens Condition $[A1']$, however 
it does not affect the proof in Yoshida \cite{Yoshida2011} because 
$S_T(r;\xi_0)=\Omega$ when $r>b_T(\xi_0)^{(1-\alpha)/2}$ in the notation of the paper. 
\end{remark}
}

{\colorr
Define $r_T(u)$ by 
\bea\label{202006010116}
r_T(u) 
&=& 
\left\{\begin{array}{ll}
\bbH_T(\theta^*+a_Tu)-\bbH_T(\theta^*)-\bigg(\Delta_T[u]-\half\Gamma[u^{\otimes2}]\bigg)
& ( u\in\bbU_T) \\
1 & (u\not\in\bbU_T) 
\end{array}\right.
\eea
}
\begin{en-text}
{\colorr 
Then $r_T(u)$ has the expression 
\bea\label{0302230320} 
r_T(u) 
&=&
\int_0^1 (1-s)\big\{\Gamma[u^{\otimes2}]
-\Gamma_T(\theta^*+sa_Tu)[u^{\otimes2}]\big\}ds
\eea
if $u\in\bbR^\sfp$ and $|a_Tu|<{\tt r}_0$. 
}
\end{en-text}
\begin{proposition}\label{202006010143}
Suppose that 
\bea\label{202006010144}
\sup_{u\in\bbU_T\cap U(0,K)}
\big|\Gamma_T(\theta^*+a_T u)-\Gamma\big| &\to^p& 0\quad(T\to\infty)
\eea
for every $K>0$. 
Then the random field $\bbZ_T$ is locally asymptotically quadratic at $\theta^*$, that is, 
\bea\label{202006010303} 
\bbZ_T(u) 
&=&
\exp\bigg(\Delta_T[u]-\half\Gamma[u^{\otimes2}]+r_T(u)\bigg)
\quad(u\in\bbU_T)
\eea
and $r_T(u)\to^p0$ as $T\to\infty$ for every $u\in\bbR^\sfp$. 
\end{proposition}
\proof
By definition of $r_T(u)$, Equation (\ref{202006010303}) holds for $u\in\bbU_T$. 
For each $u\in\bbR^\sfp$, there is a number $T_u$ such that 
$a_Tu\in U(0,{\tt r}_0)$ for all $T\geq T_u$. 
Then $r_T(u)$ admits the expression 
\bea\label{202006010301}
r_T(u) 
&=&
-\int_0^1(1-s)\big\{\Gamma_T(\theta^*+sa_Tu)-\Gamma\big\}ds[u^{\otimes2}].
\eea
{\colorr Therefore $r_T(u)\to^{{\sblue p}}0$} as $T\to\infty$ by (\ref{202006010144}). 
\qed\halflineskip
\begin{remark}\label{202006010221}\rm 
We have 
\bea\label{202006010213}
\sup_{u\in\bbU_T\cap\hspace{1pt} U(0,K)}
\big|\Gamma_T(\theta^*+a_T u)-\Gamma_T(\theta^*)\big|
&\to^p& 
0
\eea
as $T\to\infty$ under $[S3]$ (iii) since 
\beas&&
\limsup_{T\to\infty}
\left\|\sup_{u\in\bbU_T\cap\hspace{1pt} U(0,K)}
\big|\Gamma_T(\theta^*+a_T u)-\Gamma_T(\theta^*)\big|\right\|_{M_3} 
\nn\\&\leq&
\limsup_{T\to\infty}\left\|\sup_{v\in\bbR^\sfp:\>|v|\leq 1}
\big|\Gamma_T\big(\theta^*+K|a_T| v\big)-\Gamma_T(\theta^*)\big|\right\|_{M_3} 
\nn\\&\leq&
\limsup_{T\to\infty}O(K|a_T|)
\yeq 0. 
\eeas
On the other hand, 
\bea\label{202006010216}
\Gamma_T(\theta^*)-\Gamma &\to^p& 0
\eea
as $T\to\infty$ under $[S3]$ (iv). 
Therefore, 
$\bbZ_T$ is locally asymptotically quadratic at $\theta^*$ 
if $[S3]$ (iii) and (iv) are satisfied since 
the convergence (\ref{202006010144}) holds under $[S3]$ (iii) and (iv), 
though these conditions are too sufficient for (\ref{202006010144}). 
\end{remark}
\begin{en-text}
\bea\label{201912011356} 
\bbZ_T(u) &=& \exp\bigg(\Delta_T[u]-\half\Gamma[u^{\otimes2}]+r_T(u)\bigg)
\quad(u\in\bbU_T)
\eea
It is said that 
$\bbZ_T$ is locally asymptotically quadratic (LAQ) at $\theta^*$ 
if $r_T(u)\to^p0$ as $T\to\infty$ for every $u\in\bbR^\sfp$, and hence 
$\log\bbZ_T(u)$ is asymptotically approximated by a random quadratic function of $u$. 

We will confine our attention to a very standard case where 
$\bbZ_T$ is locally asymptotically mixed normal, though the general theory of the quasi-likelihood analysis 
is framed more generally.  
\end{en-text}

Let {\colorr$\Delta$} be a $\sfp$-dimensional random vector on some 
extension of $(\Omega,\calf,P)$. 
Define a random field $\bbZ$ on $\bbR^\sfp$ by 
\bea\label{202006010405} 
\bbZ(u) 
&=& 
\exp\bigg(\Delta[u]-\half\Gamma[u^{\otimes2}]\bigg)
\eea
for $u\in\bbR^\sfp$. 
Let {\colorr$\hat{u}=\Gamma^{-1}\Delta$}.

Any measurable mapping $\hat{\theta}_T^M:\Omega\to\overline{\Theta}$ is called 
a quasi-maximum likelihood estimator (QMLE) for $\bbH_T$ if 
\bea\label{0302220616}
\bbH_T(\hat{\theta}_T^M) &=& \max_{\theta\in\overline{\Theta}}\bbH_T(\theta).
\eea
Since $\bbH_T$ is continuous on the compact $\overline{\Theta}$, such a measurable function always exists, 
which is ensured by the measurable selection theorem. 
Uniqueness of $\hat{\theta}_T^M$ is not assumed. 
Let $\hat{u}_T^M=a_T^{-1}(\hat{\theta}_T^M-\theta^*)$ for the QMLE $\hat{\theta}_T^M$. 

Let $\calg$ be a $\sigma$-field such that 
$\sigma[\Gamma]\subset\calg\subset\calf$. 
{\colorr 
It is said that 
a sequence $(V_T)_{T\in\bbT}$ of random variables taking values 
in a metric space $S$ equipped with the Borel $\sigma$-field 
converges $\calg$-stably to an $S$-valued random variable $V_\infty$ defined on an extension 
of $(\Omega,\calf,P)$ 
if $(V_T,\Psi)\to^d(V_\infty,\Psi)$ as $T\to\infty$ for any $\calg$-measurable random variable $\Psi$. 
The $\calg$-stable convergence is denoted by 
$\to^{d_s(\calg)}$. 
}
\begin{theorem}\label{201911221015}
Let $L>p>0$. 
Suppose that Conditions $[S1]$, $[S2]$ and $[S3]$ are satisfied and that 
\bea\label{202006010257}
\Delta_T\to^{d_s(\calg)} \Delta
\eea
as $T\to\infty$. 
Then 
\bd
\im[(a)]  As $T\to\infty$, 
\bea\label{202006120207}
\hat{u}_T^M-\Gamma^{-1}\Delta_T  &\to^p& 0. 
\eea
\im[(b)] As $T\to\infty$, 
\bea\label{202006010421}
E\big[f(\hat{u}_T^M)\Phi\big] &\to& \bbE\big[f(\hat{u})\Phi\big]
\eea
for any bounded $\calg$-measurable random variable $\Phi$ and 
any $f\in C(\bbR^\sfp)$ satisfying\\ $\limsup_{|u|\to\infty}|u|^{-p}|f(u)|<\infty$. 
\ed
\end{theorem}
\proof 
As mentioned in Remark \ref{202006010221}, the convergence (\ref{202006010144}) holds for every $K>0$ 
under $[S3]$ (iii) and (iv). 
Then the representation (\ref{202006010301}) ensures 
\bea\label{202006010308}
\sup_{u:\>|u|\leq R}|r_T(u)| &\to^p& 0\quad(T\to\infty)
\eea
for every $R>0$. 
The space $C(\ol{U(0,R)})$ of continuous function on $\ol{U(0,R)}$ is equipped with the supremum norm. 
Combining the representation (\ref{202006010303}) of $\bbZ_T$ with the convergences 
(\ref{202006010257}) and (\ref{202006010308}), 
{\colorr by estimating the modulus of continuity of $\log\bbZ_T$ on $\ol{U(0,R)}$,}
we obtain 
tightness of the family $\big\{\bbZ_T|_{\ol{U(0,R)}}\big\}_{T\geq T_0}$ for some $T_0\in\bbT$, which 
yields the convergence 
\bea\label{202006010314}
\bbZ_T|_{\ol{U(0,R)}} &\to^d
& \bbZ|_{\ol{U(0,R)}} 
\eea
in $C(\ol{U(0,R)})$ as $T\to\infty$ for every $R>0$ 
because the finite dimensional 
convergence $\bbZ_T\to^{d_f}\bbZ$ 
is given by (\ref{202006010257}), (\ref{202006010308}) and (\ref{202006010303}). 

Let $F$ be any closed set in $\bbR^\sfp$. Then
\bea\label{202006010356} 
\limsup_{T\to\infty}P\big[\hat{u}_T^M\in F\big]
&\leq&
\limsup_{T\to\infty}P\big[\hat{u}_T^M\in F\cap\ol{U(0,R)}\>\big]+\limsup_{T\to\infty}P\big[\hat{u}_T^M\in\ol{\bbV_T(R)}\>\big]
\nn\\&\leq&
\limsup_{T\to\infty}P\bigg[\sup_{u\in F\cap\ol{U(0,R)}}\bbZ_T(u)- \sup_{u\in F^c\cap\ol{U(0,R)}}\bbZ_T(u)\geq0\bigg]
\nn\\&&
+\limsup_{T\to\infty}P\bigg[\sup_{u\in{\bbV}_T(R)}\bbZ_T(u)\geq1\bigg]
\nn\\&\leq&
P\bigg[\sup_{u\in F\cap\ol{U(0,R)}}\bbZ(u)- \sup_{u\in F^c\cap\ol{U(0,R)}}\bbZ(u)\geq0\bigg]
+\frac{C_L}{R^L}
\eea
by the convergence (\ref{202006010314}) and the polynomial type large deviation inequality (\ref{202005311807}) 
given by Theorem \ref {0302201901}. 
Let $R\to\infty$ in (\ref{202006010356}) to obtain 
\bea\label{202006010406}
\limsup_{T\to\infty}P\big[\hat{u}_T^M\in F\big]
&\leq&
P\bigg[\sup_{u\in F}\bbZ(u)- \sup_{u\in F^c}\bbZ(u)\geq0\bigg]
\yleq 
P\big[\hat{u}\in F\big]. 
\eea
Here 
{\colorr the positivity of $\Gamma$ given by $[S2]$ (ii) (Remark \ref{0302240406}) 
was used 
for the first inequality, and} 
the last inequality is by the uniqueness of the maximum point of the random field $\bbZ$ defined by (\ref{202006010405}). 
Inequality (\ref{202006010406}) shows the convergence $\hat{u}_T^M\to^d\hat{u}$ as $T\to\infty$. 
%

From the convergence of $\hat{u}_T^M$, in particular $\hat{\theta}_T^M\to^p\theta^*$, and when $\hat{\theta}_T^M\in U(\theta^*,{\tt r}_0)$, one has 
\beas 
\Delta_T
&=& 
\int_0^1{\colorr\Gamma_T}\big(\theta_T^*+s(\hat{\theta}_T^M-\theta^*)\big)ds\>\hat{u}_T^M
\eeas
{\colorr since $\partial_\theta\bbH_T(\hat{\theta}_T^M)=0$.}
Then we obtain (\ref{202006120207}) from $[S3]$ (iii) and (iv). 
{\colorr 
The $\calg$-stable convergence 
\bea\label{0302231624}
\hat{u}_T^M\to^{d_s(\calg)}\hat{u}
\eea 
follows from (\ref{202006010257}). 
}

As already used in the above argument, 
\bea
P\big[|\hat{u}_T^M|\geq r\big]
&\leq&
P\bigg[\sup_{u\in\bbV_T(r)}\bbZ_T(u)\geq1\bigg]\yleq \frac{C_L}{r^L}
\eea
for all $T\in\bbT$ and $r>0$. 
Therefore, 
\beas 
\sup_{T\in\bbT} E\big[|\hat{u}_T^M|^q\big] &<& \infty
\eeas
for any constant $q$ such that $L>q>p$. 
This means the family $\big\{f(\hat{u}_T^M)\big\}_{T>0}$ is uniformly integrable. 
Consequently, we obtain (\ref{202006010421}) 
{\colorr 
from (\ref{0302231624}). }
\qed\halflineskip
{\colorr
\begin{remark}\rm
(i) 
The convergence (\ref{202006010421}) holds for non-bounded $\Phi$ 
if $\Phi$ has the dual integrability for $f(\hat{u}^M_T)$. 
For example,
the convergence holds for $\Phi\in L^{r}(\calg)$  for some $r>1$
if $\limsup_{|u|\to\infty}|u|^{-p(r-1)/r}|f(u)|<\infty$.  
(ii) The asymptotic equivalence (\ref{202006120207}) between 
$\hat{u}^M_T$ and $\Gamma^{-1}\Delta_T$ is called 
the first-order efficiency in particular for the maximum likelihood estimator. 
This relation is useful when one considers a joint convergence of $\hat{u}^M_T$ 
with other variables. 
Such an asymptotic representation of the error is useful in analysis of a model having multi-scaled parameters.  
\end{remark}
}

The quasi-likelihood analysis enables us to derive asymptotic properties of the Bayesian estimator, 
as well as the quasi-maximum likelihood estimator. 
The mapping 
\bea\label{0302220626}
\hat{\theta}_T^B 
&=&
\bigg[\int_\Theta\exp\big(\bbH_T(\theta)\big)\varpi(\theta)d\theta\bigg]^{-1}
\int_\Theta\theta\exp\big(\bbH_T(\theta)\big)\varpi(\theta)d\theta
\eea
is called a quasi-Bayesian estimator (QBE) with respect to the prior density $\varpi$. 
The QBE $\hat{\theta}_T^B$ takes values in the convex-hull of $\Theta$. 
When the $\bbH_T$ is the log likelihood function, the QBE coincides with the Bayesian estimator with respect to 
the quadratic loss function. 
We will assume $\varpi$ is continuous and 
$0<\inf_{\theta\in\Theta}\varpi(\theta)\leq\sup_{\theta\in\Theta}\varpi(\theta)<\infty$. 
{\sblue Let $\hat{u}^B_T=a_T^{-1}\big(\hat{\theta}_T^B-\theta^*\big)$.}
%
%

\begin{theorem}\label{201912011335}
\bd
\im[{\bf (I)}] 
Let $L>1$. Suppose that Conditions $[S1]$, $[S2]$ and $[S3]$ are satisfied 
and that the convergence (\ref{202006010257}) holds as $T\to\infty$. 
Then 
\bd
\im[{\bf (a)}] As $T\to\infty$, 
\bea\label{202006120231}
\hat{u}_T^B-\Gamma^{-1}\Delta_T &\to^p& 0.
\eea

\im[{\bf (b)}]  
As $T\to\infty$,
\bea\label{202006041238}
\hat{u}_T^B &\to^{{\colorr d_s(\calg)}} & \hat{u} .
\eea

\ed

\im[{\bf (II)}]
{\colorr Let $p\geq0$ and $L>(p+1)\vee2$.} 
Suppose that Conditions $[S1]$, $[S2]$ and $[S3]$ are satisfied and that the convergence (\ref{202006010257}) holds as $T\to\infty$. 
{\colorr
Moreover, suppose that 
there exist positive constants $q$, $c_0$, $\delta$ and $T_0\in\bbT$ such that $q>\sfp$ 
and 
\bea\label{202006060415}
\sup_{T\geq T_0}E\big[\big|\bbH_T(\theta^*+a_Tu)-\bbH_T(\theta^*)\big|^q\big] 
&\leq&
c_0|u|^q
\eea
for all $u\in U(0,\delta)$.}
%
Then (\ref{202006120231}) holds, and moreover, 
{\colorr
\bea\label{202006120230}
E\big[f(\hat{u}_T^B)\Phi\big] &\to& \bbE\big[f(\hat{u})\Phi\big]
\eea
as $T\to\infty$ 
for any $\calg$-measurable bounded random variable $\Phi$ and}
any $f\in C(\bbR^\sfp)$ satisfying\\ $\limsup_{|u|\to\infty}|u|^{-p}|f(u)|<\infty$. 
\ed
\end{theorem}
\halflineskip
\begin{remark}\label{202108251451}
\rm 
{\colorr
In Theorem \ref{201912011335}, 
we implicitly assume that $T_0$ is sufficiently large so that 
$U(0,\delta)\subset\bbU_T$ for all $T\geq T_0$, 
and the left-hand side of (\ref{202006060415}) makes sense.} 
\end{remark}
\begin{remark}\label{202006060419}\rm 
Condition (\ref{202006060415}) holds
{\colorr under any one of the following conditions: 
\bd
\im[(i)] There exist constants $q>\sfp$, $\delta>0$ and $T_0\in\bbT$ such that 
{\sblue 
\beas 
\sup_{T\geq T_0}\sup_{u\in \bbU_T\cap U(0,\delta)}\|\Gamma_T(\theta^*+a_Tu)\|_q
&<&
\infty
\eeas
}\noindent
and $\sup_{T\geq T_0}\|\Delta_T\|_q<\infty$. 
\im[(ii)]$|\Gamma|\in L^q$, and $M_1$, $M_3,$ and $M_4$ appearing in $[S3]$ 
satisfy 
$
M_1\wedge M_3\wedge M_4\ygeq q
$.
\ed
}\noindent
This follows from the representation (\ref{202006010303}) of $\bbZ_T(u)$ and 
the representation (\ref{202006010301}) of $r_T(u)$. 
\end{remark}
\halflineskip\noindent
{\it Proof of Theorem \ref{201912011335}.} 
%
(I) 
We obtain a polynomial type large deviation inequality from 
Theorem \ref{0302201901}: for any $D>0$, 
there exist positive constants 
$C_1$ and $C_2$ such that 
\bea\label{201912011344} 
P\bigg[\sup_{u\in\bbV_T(r)}\bbZ_T(u)\geq C_1r^{-D}\bigg]&\leq& C_2r^{-L}
\eea
for all $T$ and $r>0$. 
{\colorr Choose a number $D$ such that $D>\sfp+(p\vee1)$.} 
 If we take a sufficiently large constant $C_1'$, then 
\beas &&
P\bigg[\int_{{\colorr\bbV_T(r)}}
(1+|u|)\bbZ_T(u)du
> C_1' \sum_{\ell=0}^\infty(r+\ell)^{\sfp-D}\bigg]
\\&\leq& 
\sum_{\ell=0}^\infty P\bigg[
\int_{\{r+\ell\leq|u|< (r+\ell+1)\}{\colorr\cap\bbU_T}}(1+|u|)\bbZ_T(u)du
> C_1' (r+\ell)^{\sfp-D}\bigg]
\\&\leq& 
\sum_{\ell=0}^\infty P\bigg[
\sup_{\{r+\ell\leq|u|< (r+\ell+1)\}{\colorr\cap\bbU_T}}\bbZ_T(u)
> C_1(r+\ell)^{-D}\bigg]
\\&\leq& 
C_2\sum_{\ell=0}^\infty(r+\ell)^{-L}
\eeas
for $T\in\bbT$ and $r>1$, by (\ref{201912011344}). 
{\colorr Since $D-\sfp>1$ and $L>1$ by assumption,}
there exist positive constants $C_3$ and $\ep$ (independent of $(r,T)$) such that 
\bea\label{202006041205} 
P\bigg[\int_{\bbV_T(r)}(1+|u|)\bbZ_T(u)du>C_3r^{-\ep}\bigg]
&\leq& 
C_3r^{-\ep}\quad(r>0,\>T\in\bbT)
\eea

\begin{en-text}
{\color{gray}
Let $R>r>0$. For large $T$, $\bbZ_T(u)$ is well defined on $\{u;\>r\leq |u|\leq R\}$, and 
the continuous mapping theorem applied to 
the convergence (\ref{202006010314}), that was shown in the proof of Theorem \ref{201911221015}, 
gives the inequality 
\beas 
P\bigg[\int_{r\leq|u|}(1+|u|)\bbZ(u)du>C_3r^{-\ep}\bigg]
&\leq&
\lim_{R\to\infty}
P\bigg[\int_{r\leq|u|\leq R}(1+|u|)\bbZ(u)du>C_3r^{-\ep}\bigg]
\nn\\&\leq&
\lim_{R\to\infty}\liminf_{T\to\infty}P\bigg[\int_{r\leq|u|\leq R}(1+|u|)\bbZ_T(u)du>C_3r^{-\ep}\bigg]
\nn\\&\leq&
C_3r^{-\ep}\quad(r>0)
\eeas
from (\ref{202006041205}). 
Then the Borel-Cantelli lemma gives 
\beas 
\int_{\bbR^\sfp}(1+|u|)\bbZ(u)du &<& \infty \quad a.s. 
\eeas
}
\end{en-text}

The variable $\hat{u}_T^B$ has the expression 
\bea\label{202006061303} 
\hat{u}_T^B
&=& 
\bigg(\int_{\bbU_T}\bbZ_T(u)\varpi(\theta^*+a_Tu)du\bigg)^{-1}
\int_{\bbU_T} u\bbZ_T(u)\varpi(\theta^*+a_Tu)du.
\eea
{\colorr For $g(u)=(1,u)$,} 
let 
\beas &&
X_T=\int_{\bbU_T}g(u)\bbZ_T(u)\varpi(\theta^\dagger_T(u))du, \quad 
X_{T,r}=\int_{\bbU_T\cap U(0,r)}g(u)\bbZ_T(u)\varpi(\theta^\dagger_T(u))du, \quad 
\nn\\&&
W_{T,r}=\int_{\bbV_T(r)}g(u)\bbZ_T(u)\varpi(\theta^\dagger_T(u))du,\quad
\nn\\&&
X_\infty=\int_{\bbR^\sfp}g(u)\bbZ(u)\varpi(\theta^*)du, \quad 
X_{\infty,r}=\int_{U(0,r)}g(u)\bbZ(u)\varpi(\theta^*)du, \quad 
\eeas
where $\theta^\dagger_T(u)=\theta^*+a_Tu$ and $\bbZ$ is given by (\ref{202006010405}). 
Then $X_T=X_{T,r}+W_{T,r}$ and 
{\colorr the following properties hold.}
\bi
\im[(i)] For any $\eta>0$, there exists $r_0>0$ such that 
$\sup_{T\in\bbT}P[|W_{T,r}|>\eta]<\eta$ for all $r\geq r_0$. 
\im[(ii)] For every $r>1$, $X_{T,r}\to^d X_{\infty,r}$ as $T\to\infty$. 
\im[(iii)] $X_{\infty,r}\to^d  X_\infty$ as $r\to\infty$. 
\ei
Indeed, (i) follows from (\ref{202006041205}), (ii) from the convergence (\ref{202006010314}), and 
(iii) is obvious. 
{\colorr Therefore 
\bea\label{0302231243}
X_T &\to^d& X_\infty
\eea
as $T\to\infty$. }

Denote $X_T=(X_T^{(0)},X_T^{(1)})$ and $X_{T,r}=(X_{T,r}^{(0)},X_{T,r}^{(1)})$. 
{\colorr We will consider sufficiently large $T$ such that $\bbU_T\supset U(0,1)$. 
}
Let 
\beas 
A_T &=& \big(|X^{(1)}_T|+X^{(0)}_T)\bigg(X^{(0)}_T\int_{U(0,1)}\bbZ_T(u)du\>\inf_{\theta\in\Theta}\varpi(\theta)\bigg)^{-1}. 
\eeas
Then 
\beas
A_T&\geq & \frac{|X^{(1)}_T|+X^{(0)}_T}{X^{(0)}_TX^{(0)}_{T,r}}
\eeas
for all $r\geq1$. 
We have 
\beas 
\left| \big(X_T^{(0)}\big)^{-1}X_T^{(1)}-\big(X_{T,r}^{(0)}\big)^{-1}X_{T,r}^{(1)}\right|
&\leq&
|W_{T,r}|A_T
\eeas
for all $r\geq1$. 
Let $\ep>0$. Then there exists a positive number $\eta>0$ such that 
\beas 
\limsup_{T\to\infty}P\bigg[A_T>\frac{1}{\eta}\bigg] &<& \frac{\ep}{4}
\eeas
since the family $\{A_T\}_{T\geq T_1}$ is tight for some $T_1\in\bbT$ 
by (\ref{0302231243}) 
and (\ref{202006010314}). 
For the pair $(\ep,\eta)$, there exists $r_0=r_0(\ep,\eta)\geq1$ such that 
\beas 
\limsup_{T\to\infty}P\bigg[|W_{T,r}|>\frac{\ep\eta}{4}\bigg] &<& \frac{\ep}{4}\quad(r\geq r_0)
\eeas
{\colorr by the property (i) mentioned just before (\ref{0302231243}). 
In what follows, we fix an $r\geq r_0$. }
Then
\bea\label{202006120456} &&
\limsup_{T\to\infty}
P\bigg[\left| \big(X_T^{(0)}\big)^{-1}X_T^{(1)}-\big(X_{T,r}^{(0)}\big)^{-1}X_{T,r}^{(1)}\right|>\frac{\ep}{4}\bigg]
\nn\\&\leq&
\limsup_{T\to\infty}P\bigg[|W_{T,r}|>\frac{\ep\eta}{4}\bigg]
+\limsup_{T\to\infty}P\bigg[A_T>\frac{1}{\eta}\bigg]
\><\>\frac{\ep}{2}.
\eea

We have 
\beas &&
\bigg|\int_{U(0,r)}u^i\bbZ_T(u)\varpi(\theta_T^\dagger(u))du-\int_{U(0,r)}u^i\bbZ_T(u)\varpi(\theta^*)du\bigg|
\nn\\&\leq&
\int_{U(0,r)}(1+|u|)\bbZ_T(u)du\>
\nn\\&&\times
{\colorr
\sup\bigg\{|\varpi(\theta)-\varpi(\theta^*)|;\>\theta\in\Theta,
\>|\theta-\theta^*|\leq|a_T|r\bigg\}
}
\eeas
where 
$u^0=1$ and $u^1=u$. 
Therefore 
\bea\label{202006120515} &&
\bigg(\int_{U(0,r)}\bbZ_T(u)\varpi(\theta_T^\dagger(u))du\bigg)^{-1}\int_{U(0,r)}u\bbZ_T(u)\varpi(\theta_T^\dagger(u))du
\nn\\&&\hspace{30pt}
-
\bigg(\int_{U(0,r)}\bbZ_T(u)du\bigg)^{-1}\int_{U(0,r)}u\bbZ_T(u)du
\>\to^p\>
0
\eea
as $T\to\infty$. 

Moreover, we have 
\beas &&
\bigg|\int_{U(0,r)}u^i\bbZ_T(u)du-\int_{U(0,r)}u^i\hat{\bbZ}_T(u) du\bigg|
\nn\\&=&
\bigg|\int_{U(0,r)}u^i\exp\big(\Delta_T[u]-2^{-1}\Gamma[u^{\otimes2}]+r_T(u)\big) du
-\int_{U(0,r)}u^i\exp\big(\Delta_T[u]-2^{-1}\Gamma[u^{\otimes2}]\big) du\bigg|
\nn\\&\leq& 
\int_{U(0,r)}(1+|u|)\hat{\bbZ}_T(u) du
\>\sup_{u\in U(0,r)}\big|e^{r_T(u)}-1\big|
\eeas
where 
\beas 
\hat{\bbZ}_T(u)
&=&
\exp\big(\Delta_T[u]-2^{-1}\Gamma[u^{\otimes2}]\big).
\eeas
Therefore
\bea\label{202006120516}
\bigg(\int_{U(0,r)}\bbZ_T(u)du\bigg)^{-1}\int_{U(0,r)}u\bbZ_T(u)du
-\bigg(\int_{U(0,r)}\hat{\bbZ}_T(u)du\bigg)^{-1}\int_{U(0,r)}u\hat{\bbZ}_T(u)du
&\to^p&
0
\nn\\&&
\eea
as $T\to\infty$ 
{\colorr 
thanks to the convergence (\ref{202006010308}).}

Let 
\beas &&
\hat{X}_T=(\hat{X}_T^{(0)},\hat{X}_T^{(1)})=\int_{\bbR^\sfp}g(u)\hat{\bbZ}_T(u)du, \quad 
\hat{X}_{T,r}=(\hat{X}_{T,r}^{(0)},\hat{X}_{T,r}^{(1)})=\int_{U(0,r)}g(u)\hat{\bbZ}_T(u)du, \quad 
\nn\\&&
\hat{W}_{T,r}=\int_{\bbR^\sfp\setminus U(0,r)}g(u)\hat{\bbZ}_T(u)du.
\eeas
Then
\beas 
\left| \big(\hat{X}_T^{(0)}\big)^{-1}\hat{X}_T^{(1)}-\big(\hat{X}_{T,r}^{(0)}\big)^{-1}\hat{X}_{T,r}^{(1)}\right|
&\leq&
|\hat{W}_{T,r}|\hat{A}_T
\eeas
for all $r\geq1$, where 
\beas 
\hat{A}_T &=& \big(|\hat{X}^{(1)}_T|+\hat{X}^{(0)}_T\big)\bigg(\hat{X}^{(0)}_T\int_{U(0,1)}\hat{\bbZ}_T(u)du\>\inf_{\theta\in\Theta}\varpi(\theta)\bigg)^{-1}
\nn\\&\geq&
\frac{|\hat{X}^{(1)}_T|+\hat{X}^{(0)}_T}{\hat{X}^{(0)}_T\hat{X}^{(0)}_{T,r}}
\eeas
{\colorr 
Positive-definiteness of $\Gamma$ (Remark \ref{0302240406} (iii)) 
and the tightness due to the convergence of $\Delta_T$ 
show that for any $\eta>0$, there exist $r_1>0$ and $T_2\in\bbT$ such that 
\bea\label{0302240419}
\sup_{T\geq T_2}P[\hat{W}_{T,r}>\eta]<\eta\quad(r\geq r_1).
\eea 
}
In the same way as we showed (\ref{202006120456}), 
\bea\label{202006120457} 
\limsup_{T\to\infty}
P\bigg[\left| \big(\hat{X}_T^{(0)}\big)^{-1}\hat{X}_T^{(1)}-\big(\hat{X}_{T,r}^{(0)}\big)^{-1}\hat{X}_{T,r}^{(1)}\right|>\frac{\ep}{4}\bigg]
&<&
\frac{\ep}{2}
\eea
{\sred for some $r$.}
{\colorr
To obtain (\ref{202006120457}) 
by using  (\ref{0302240419}) and the tightness of $\{\hat{A}_T\}_{T\geq T_3}$ 
for some $T_3\in\bbT$,  
we replace $r$ by a larger number, if necessary.}

Combining 
(\ref{202006120456}), (\ref{202006120515}), (\ref{202006120516}) and (\ref{202006120457}), 
we obtain 
\bea\label{202006120506} 
\limsup_{T\to\infty}
P\bigg[\left| \big(X_T^{(0)}\big)^{-1}X_T^{(1)}-\big(\hat{X}_T^{(0)}\big)^{-1}\hat{X}_T^{(1)}\right|>\ep\bigg]
&<&
\ep.
\eea
This completes the proof of (\ref{202006120231}) 
since $\hat{u}_T^B=\big(X_T^{(0)}\big)^{-1}X_T^{(1)}$ and 
$\big(\hat{X}_T^{(0)}\big)^{-1}\hat{X}_T^{(1)}=\Gamma^{-1}\Delta_T$. 
From (\ref{202006120231}), we obtain 
(\ref{202006041238}). %
\begin{en-text}
Therefore, 
\bea\label{0302231243}
X_T&\to^{{\colorr d_s(\calg)}} &X_\infty
\eea
as $T\to\infty$, in particular, 
\bea\label{202006041238}
\hat{u}_T^B &\to^{{\colorr d_s(\calg)}} & \hat{u} 
\eea
as $T\to\infty$ thanks to (\ref{202006061303}) 
since 
\beas 
\hat{u}\yeq \Gamma^{-1}\Delta
&=&
\bigg(\int_{\bbR^\sfp}\bbZ(u)du\bigg)^{-1}
\int_{\bbR^\sfp} u\bbZ(u)du.
\eeas
\end{en-text}

\halflineskip
\noindent (II) 
There exists a number $p_*$ such that 
\beas 
p_*\geq1,\quad
p<p_*<(D-\sfp)\wedge(L-1)
\eeas
since {\colorr $p\geq0$ and $L>(p+1)\vee2$}. 
Then the following estimates are standard: 
\beas 
E\big[|\hat{u}_T^B|^{p_*}]
&\leq&
E\bigg[\bigg(\int_{\bbU_T}\bbZ_T(u)\varpi(\theta^\dagger_T(u))du\bigg)^{-1}
\int_{\bbU_T}|u|^{p_*}\bbZ_T(u)\varpi(\theta^\dagger_T(u))du\bigg]
\\&\leq&
C(\varpi)\sum_{r=0}^\infty
(r+1)^{p_*} E\bigg[\bigg(\int_{\bbU_T}\bbZ_T(u)du\bigg)^{-1}
\int_{\{u;r<|u|\leq r+1\}\cap\bbU_T}\bbZ_T(u)du\bigg]
\\&\leq&
C(\varpi)\big(1+\Phi_{1,T}+\Phi_{2,T}\big)
\eeas
for sme constant $C(\varpi)$, 
where 
\beas 
\Phi_{1,T}
&=&
\sum_{r=1}^\infty(r+1)^{p_*}
E\bigg[\bigg(\int_{\bbU_T}\bbZ_T(u)du\bigg)^{-1}
\int_{\{u;r<|u|\leq r+1\}\cap\bbU_T}\bbZ_T(u)du 
\nn\\&&\hspace{80pt}\times
\>1_{\bigg\{\int_{\{u;r<|u|\leq r+1\}\cap\bbU_T}
\bbZ_T(u)du>\frac{C_1'}{r^{D-\sfp+1}}
\bigg\}}\bigg]
\eeas
and 
\beas 
\Phi_{2,T}
&=&
\sum_{r=1}^\infty(r+1)^{p_*} 
E\bigg[\bigg(\int_{\bbU_T}\bbZ_T(u)du\bigg)^{-1}
\int_{\{u;r<|u|\leq r+1\}\cap\bbU_T}\bbZ_T(u)du 
\nn\\&&\hspace{80pt}\times
\>1_{\bigg\{\int_{\{u;r<|u|\leq r+1\}\cap\bbU_T}
\bbZ_T(u)du\leq\frac{C_1'}{r^{D-\sfp+1}}
\bigg\}}\bigg].
\eeas
Take a sufficiently large number $C_1'$. 
Since the integrand of the expectation of $\Phi_{1,T}$ is not greater than one, we obtain  
\beas 
\Phi_{1,T}
&\leq&
\sum_{r=1}^\infty(r+1)^{p_*}
P\bigg[\int_{\{u;r<|u|\leq r+1\}\cap\bbU_T}
\bbZ_T(u)du>\frac{C_1'}{r^{D-\sfp+1}}\bigg]
\\&\simleq&
\sum_{r=1}^\infty r^{-(L-{p_*})}
\eeas
thanks to the polynomial type large deviation inequality (\ref{201912011344}). 
For $\Phi_{2,T}$, 
\beas 
\Phi_{2,T}
&\simleq&
\sum_{r=1}^\infty r^{-(D-{p_*}-\sfp+1)}E\bigg[\bigg(\int_{\bbU_T}\bbZ_T(u)du\bigg)^{-1}\bigg].
\eeas
Therefore,  the family $\{|\hat{u}_T^B|^p\}_{{\sred T\geq T_0}}$ is uniformly integrable if 
\bea\label{202006060320}
\sup_{T{\colorr \geq T_0}}E\bigg[\bigg(\int_{U(0,\delta)}\bbZ_T(u)du\bigg)^{-1}\bigg]
&<& 
\infty
\eea
since {\colorr$U(0,\delta)\subset{\colorr\cap_{T\geq T_0}\bbU_T}$ 
(see Remark \ref{202108251451}) and then} 
\beas 
\sup_{\sred T\geq T_0}E\big[|\hat{u}_T^B|^{p_*}]
&<& 
\infty. 
\eeas
[We note that the family $\{\hat{u}_T^B\}_{T\in\bbT,T< T_0}$ 
is {\sred not necessarily} uniformly bounded {\sred when $\#\{T<T_0\}=\infty$} because 
{\sred $\inf_{T<T_0}|a_T|$ may be zero for example, though} 
$\hat{\theta}_T^B$ takes values 
in the convex hull of $\Theta$, that is bounded.]  
We obtain the inequality (\ref{202006060320}) 
from (\ref{202006060415}) by applying e.g. Lemma 2 of \cite{Yoshida2011}. 
The convergence (\ref{202006120230}) is now obtained from (\ref{202006041238}) and 
the uniform integrability of the family $\{|\hat{u}_T^B|^p\}_{{\sred T\geq T_0}}$. 
\qed\halflineskip
{\sred 
\begin{remark}\rm
From (\ref{202006120207}) and (\ref{202006120231}), 
the QMLE and the QBE are asymptotically equivalent in that 
$\hat{u}_T^M-\hat{u}_T^B\to^p0$ as $T\to\infty$. 
\end{remark}
\halflineskip
}

The following conditions $[T1]$ and $[T2]$ strengthen $[S2]$ and $[S3]$, respectively. 
\bd
\im[{\bf [T1]}] 
\bd
\im[{\bf (i)}] 
There exists a positive random variable $\chi_0$ and the following conditions are fulfilled. 
\bd
\im[{\bf (i-1)}] 
$
\bbY(\theta) \yeq \bbY(\theta) -\bbY(\theta^*)
\yleq -\chi_0\big|\theta-\theta^*\big|^2
$
for all $\theta\in\Theta$. 
\im[{\bf (i-2)}] 
For every $L>0$, there exists a constant $C$ such that 
\beas 
P\big[\chi_0\leq r^{-1}\big]
&\leq& \frac{C}{r^L}\quad(r>0).
\eeas
\ed
\im[{\bf (ii)}] 
For every $L>0$, there exists a constant $C$ such that 
\beas 
P\big[\lambda_{min}(\Gamma)<
r^{-1}\big]\leq\frac{C}{r^L}\quad(r>0)
\eeas
\ed
\ed
{\colorr
\begin{remark}\rm
$\chi_0^{-1}\in L^\inftym=\cap_{p>1}L^p$ under $[T1]$ (i-2). 
$|\Gamma^{-1}|\in L^\inftym$ under $[T1]$ (ii) 
since $\big(\lambda_{min}(\Gamma)\big)^{-1}\in L^\inftym$ and 
$|\Gamma^{-1}|\leq C_\sfp \big(\lambda_{min}(\Gamma)\big)^{-1}$ 
for a constant only depending on $\sfp$. 
The $L^q$-integrability of $\Gamma$ will be assumed 
when we verify (\ref{202006060415}). 
\end{remark}
}


\bd
\im[{\bf [T2]}] 
There exist positive numbers $\ep_1$ and $\ep_2$ such that 
the following conditions are satisfied for all $p>1$:
\bd
\im[{\bf (i)}] 
$
\sup_{T\in\bbT}\big\||\Delta_T|\big\|_p<\infty.
$
\im[{\bf (ii)}]  
$\ds
\sup_{T\in\bbT}
\bigg\|\sup_{\theta\in{\colorr\Theta}}
{b_T}^{{\colorr\hspace{-3pt}\ep_1}}\big|\bbY_T(\theta)-\bbY(\theta)\big|\bigg\|_p
<\infty.
$
\im[(iii)] 
$\ds 
\sup_{T\in\bbT}\bigg\|\sup_{u\in{\sred\delta^{-1}(\Theta-\theta^*),}\>|u|\leq1}
\big|\Gamma_T(\theta^*+\delta u)-\Gamma_T(\theta^*)\big|\bigg\|_p = O(\delta)\quad(\delta\down0). 
$
\im[{\bf (iv)}] 
$\ds
\sup_{T\in\bbT}\big\|
{b_T}^{{\colorr\hspace{-3pt}\ep_2}}\big|\Gamma_T{\colorr(\theta^*)}-\Gamma\big|\big\|_p
< \infty.
$
\ed
\ed

Condition $[T2]$ requires any order of moments of variables. 
In applications of the QLA to inference for stochastic differential equations, usually one needs not be nervous about the existence of moments of arbitrary order. 
Denote by $f\in C_p(\bbR^\sfp)$ the set of continuous functions 
of at most polynomial growth. 
We can further simplify Theorems \ref{201911221015} and \ref{201912011335} 
under Conditions $[T1]$ and $[T2]$. 
\begin{theorem}\label{0302220658a}
Suppose that Conditions $[T1]$ and $[T2]$ are satisfied and that 
the convergence (\ref{202006010257}) holds as $T\to\infty$. 
Then
\bd
\im[(a)] As $T\to\infty$, 
\bea\label{0302220728}
\hat{u}_T^M-\Gamma^{-1}\Delta_T  &\to^p& 0.
\eea
%
\im[(b)]  As $T\to\infty$, 
\bea\label{0302220656}
E\big[f(\hat{u}_T^M)\Phi\big] &\to& \bbE\big[f(\hat{u})\Phi\big]
\eea
for any $f\in C_p(\bbR^\sfp)$ and 
any $\calg$-measurable random variable $\Phi\in\cup_{p>1}L^p$. 
\ed
\end{theorem}
\proof 
There exist values of the parameters $\alpha,\beta_1,\beta_2,\rho_1$ and $\rho_2$ 
satisfying $\beta_1\in(0,\min\{\ep_2,1/2\})$, 
$1/2-\beta_2\leq\ep_1$ and Condition $[S1]$. 
Then $[S2]$ is verified for any given $L>0$ by $[T1]$. 
Condition $[T2]$ is sufficient for $[S3]$ for any $L>0$. 
Therefore we can apply 
Theorems \ref{201911221015}. 
This concludes the proof. 
\qed 

\begin{theorem}\label{0302220658}
Suppose that Conditions $[T1]$ and $[T2]$ are satisfied and that 
the convergence (\ref{202006010257}) holds as $T\to\infty$. 
{\colorr Moreover, suppose that $|\Gamma|\in L^q$ for some $q>\sfp$.} 
Then
\bd
\im[(a)] As $T\to\infty$, 
\bea\label{0302220728}
\hat{u}_T^B-\Gamma^{-1}\Delta_T  &\to^p& 0.
\eea
%
\im[(b)] As $T\to\infty$, 
\bea\label{0302220656}
E\big[f(\hat{u}_T^B)\Phi\big] &\to& \bbE\big[f(\hat{u})\Phi\big]
\eea
for any $f\in C_p(\bbR^\sfp)$ and 
any $\calg$-measurable random variable $\Phi\in\cup_{p>1}L^p$. 
\ed
\end{theorem}
\proof {\colorr We apply Theorem \ref{201912011335}. 
In particular, (\ref{202006060415}) holds now according to Remark \ref{202006060419}. 
}
\qed\halflineskip

\section{Further simplification in ergodic statistics
}
When the limit $\bbY$ of $\bbY_T$ is deterministic, 
more simplification of the theory is possible. 
In this section, we 
suppose that the random field $\bbY$ and 
the $\sfp\times\sfp$ positive-definite symmetric matrix $\Gamma$ 
are deterministic. 
Let $L$ be a positive number. 
We will consider the following conditions. 
Some simplification has been made in 
Condition $[U2]$, that is slightly different from $[S1]$ plus $[S3]$. 
\bd
\begin{en-text}
{\color{gray}
\im[{\bf[$E2$]}] 
For some constant $C_L$, 
\beas 
\sup_{T>T_0}P\big[S_T'(r)^c\big] &\leq& \frac{C_L}{r^L}\quad(r>0).
\eeas
}
\end{en-text}

\im[{\bf[U1]}] 
$\Gamma$ is positive-definite, in addition, 
there is a positive number $\chi_0$ such that 
\beas 
\bbY(\theta)\>=\>\bbY(\theta)-\bbY(\theta^*)
&\leq& -\chi_0|\theta-\theta^*|^2
\eeas
for all $\theta\in\Theta$.  
\begin{en-text}
\im[{\bf[$L3$]}] 
There exists a $C_L$ such that
\beas 
P\big[\chi_0\leq r^{-(\rho_2-\alpha\rho)}\big] &\leq& \frac{C_L}{r^L}\quad(r>0)
\eeas
and 
\beas 
P\big[\lambda_{\text{min}}(\Gamma)<4r^{-\rho_1}\big]&\leq&\frac{C_L}{r^L}\quad(r>0).
\eeas
\end{en-text}

\im[{\bf[U2]}] 
The numbers $\alpha$, $\beta_1$, 
$\beta_2$ 
and $\rho_2$ satisfy the inequalities
\bea\label{202006010103} 
&&
0<2\alpha<\rho_2,\quad
\beta_2\geq0,\quad
1-2\beta_2-\rho_2>0, \quad
0<\beta_1<1/2, 
\eea
and the following conditions are fulfilled.
\bd
\im[(i)] 
For some $M_1>L$, 
$\displaystyle
\sup_{{\sblue T\in\bbT}}\big\|\Delta_T\big\|_{M_1} <\infty.
$
\im[(ii)] For $M_2=L(1-2\beta_2-\rho_2)^{-1}$, 
\beas 
\sup_{{\sblue T\in\bbT}}\bigg\|
\sup_{{\sblue\theta\in\Theta\setminus U(\theta^*,b_T^{-\alpha/2})}}
b_T^{\half-\beta_2}\big|\bbY_T({\colorr\theta})-\bbY({\colorr\theta})\big|\bigg\|_{M_2}
&<&\infty. 
\eeas
\im[(iii)] For some $M_3>L\beta^{-1}$, 
\beas 
\sup_{{\sblue T\in\bbT}}\bigg\|\sup_{u\in{\sred\delta^{-1}(\Theta-\theta^*),}\>|u|\leq1}
\big|\Gamma_T(\theta^*+\delta u)-\Gamma_T(\theta^*)\big|\bigg\|_{M_3} &=& O(\delta)\quad(\delta\down0). 
\eeas
\im[(iv)] 
For some $M_4>L\big(2\beta_1)^{-1}(1-\alpha)$, 
\beas 
\sup_{{\sblue T\in\bbT}}\big\|b_T^{\beta_1}
\big|\Gamma_T(\theta^*)-\Gamma\big|\big\|_{M_4}
&<& \infty.
\eeas
\ed
\ed
\halflineskip
\begin{en-text}
Let $\bbU_T=\{u\in\bbR^\sfp;\>\theta^*+a_Tu\in\Theta\}$ and $\bbV_T(r)=\{u\in\bbU_T;\>|u|\geq r\}$ for $r>0$. 
Define the random field $\bbZ_T$ on $\bbU_T$ by 
\beas 
\bbZ_T(u) 
&=& 
\exp\big(\bbH_T(\theta^*+a_Tu)-\bbH_T(\theta^*)\big)
\eeas
for $u\in\bbU_T$. 
Following \cite{Yoshida2011}, 
we will give a simplified version of the polynomial type large deviation inequality for 
the random field $\bbZ_T$. 
\end{en-text}
\begin{remark}\rm 
Condition $[U1]$ is almost trivial because the function $\bbY$ should be of $C^2$ on $\Theta$ and continuous 
on the compact set $\ol{\Theta}$ and then local non-degeneracy of the information implies 
the global identifiability. 
\end{remark}
\begin{theorem}\label{201911220753}
Suppose that Conditions $[U1]$ and $[U2]$ are fulfilled for a positive constant $L$. Then 
there exists a constant $C_L$ such that 
\beas
P\bigg[\sup_{u\in\ol{\bbV}_T(r)}\bbZ_T(u)\geq \exp\big(-2^{-1}r^{2-\rho_2}\big)
\bigg]&\leq&\frac{C_L}{r^L}
\eeas
for all $T\in\bbT$ and $r>0$. 
Here the supremum on the empty set should read $-\infty$ by convention. 
\end{theorem}
\proof 
Choose a positive constant $\rho_1$ such that 
\bea\label{202005311812}&&
{\colorr0<\rho_1<\min\big\{1,\alpha/(1-\alpha),2\beta_1/(1-\alpha)\big\},\quad}
\rho_1 \leq \rho_2,\quad
\nn\\&&
M_1\ygeq M_1'\>:=\>L(1-\rho_1)^{-1},\quad
M_3\ygeq M_3'\>:=\>L(\beta-\rho_1)^{-1},\quad
\nn\\&&
M_4\ygeq M_4'\>:=\>L\big(2\beta_1(1-\alpha)^{-1}-\rho_1\big)^{-1}
\eea
for $M_1,M_2,M_3$ given in Condition $[U2]$. 
Such a positive number $\rho_1$ exists. 
It is sufficient to verify the conditions of Theorem \ref{0302201901}. 
Condition $[S1]$ is fulfilled by $[U2]$ and 
a choice of $\rho_1$ in (\ref{202005311812}). 
Condition $[S2]$ (i-1) is satisfied with $[U1]$, and 
Conditions (i-2) and (ii) of $[S2]$ are trivial because 
$\chi_0$ is a deterministic positive number and 
$\Gamma$ is positive-definite, deterministic in the present situation, respectively. 
Conditions (i)-(iv) of $[S3]$ are verified by $[U2]$ 
with $(M_1',M_2,M_3',M_4')$ for $(M_1,M_2,M_3,M_4)$ in $[S3]$. 
\qed\halflineskip

As before, the quasi-maximum likelihood estimator (QMLE) for $\bbH_T$ 
is characterized by (\ref{0302220616}). 
Theorem \ref{201911221015} is rephrased as follows with the trivial $\sigma$-field for $\calg$. 
\begin{theorem}\label{0302220650}
Let $L>p>0$. 
Suppose that Conditions $[U1]$ and $[U2]$ are satisfied and that 
\bea\label{0302220630}
\Delta_T\to^d \Delta
\eea
as $T\to\infty$. 
Then, 
\bd
\im[(a)] 
$\ds 
\hat{u}_T^M-\Gamma^{-1}\Delta_T  \to^p 0 
$
as $T\to\infty$. 
\im[(b)] 
$\ds 
E\big[f(\hat{u}_T^M)\big] \to \bbE\big[f(\hat{u})\big]
$
as $T\to\infty$ 
for 
any $f\in C(\bbR^\sfp)$ satisfying $\limsup_{|u|\to\infty}|u|^{-p}|f(u)|<\infty$. 
\ed
\end{theorem}
\proof 
Take $\rho_1$ as (\ref{202005311812}) and apply Theorem \ref{201911221015}. 
\qed\halflineskip

Consider the quasi-Beyesian estimator $\hat{\theta}_T^B$ defined by (\ref{0302220626}). 
We can rephrase Theorem \ref{201912011335} as follows. 
\begin{theorem}\label{0302220629}
{\colorr Let $p\geq0$ and $L>(p+1)\vee2$.} 
Suppose that Conditions $[U 1]$ and $[U2]$ are satisfied and that the convergence (\ref{0302220630}) holds as $T\to\infty$. 
{\colorr
Moreover, suppose that there exist positive constants $q$, $\delta$, $T_0\in\bbT$ 
and $c_0$ such that 
$q>\sfp$ and 
the inequality (\ref{202006060415}) holds for all $u\in U(0,\delta)$.}
%
Then 
\bd
\im[(a)]  
$\ds 
\hat{u}_T^B-\Gamma^{-1}\Delta_T \to^p 0
$
as $T\to\infty$. 
\im[(b)]  
{\colorr
$\ds 
E\big[f(\hat{u}_T^B)\big] \to \bbE\big[f(\hat{u})\big]
$
as $T\to\infty$ 
for any $f\in C(\bbR^\sfp)$ satisfying $\limsup_{|u|\to\infty}|u|^{-p}|f(u)|<\infty$.}
\ed
\end{theorem}
\halflineskip

As a corollary of Theorem \ref{0302220658a} (or Theorem \ref{0302220650}) and Theorem \ref{0302220658} 
(or Theorem \ref{0302220629}), 
we obtain the following result. 
\begin{theorem}\label{0302220806}
Suppose that Conditions $[U1]$ and $[T2]$ are satisfied and that 
the convergence (\ref{0302220630}) holds as $T\to\infty$. 
Then
\bd
\im[(a)] 
$\ds 
\hat{u}_T^{\sf A}-\Gamma^{-1}\Delta_T  \to^p 0
$
as $T\to\infty$ for ${\sf A}\in\{M,B\}$. 
\im[(b)]  
$\ds 
E\big[f(\hat{u}_T^{\sf A})\big] \to \bbE\big[f(\hat{u})\big]
$
as $T\to\infty$ 
for ${\sf A}\in\{M,B\}$ and any $f\in C_p(\bbR^\sfp)$. 
\ed
\end{theorem}
%
%
\bibliographystyle{spmpsci}      
\bibliography{bibtex-20200418-20200531-20201001+++}   

\end{document}